\documentclass[reqno]{amsart}
\usepackage{amscd}
\usepackage[dvips]{graphics}

\def\beq{\begin{equation}}
\def\eeq{\end{equation}}
\def\ba{\begin{array}}
\def\ea{\end{array}}
\def\S{\mathbb S}
\def\R{\mathbb R}
\def\C{\mathbb C}

\def\N{\mathbb N}

\def \ds{\displaystyle}
\def \vs{\vspace*{0.1cm}}

\newtheorem{thm}{Theorem}[section]
\newtheorem{lm}[thm]{Lemma}
\newtheorem{prop}[thm]{Proposition}

\theoremstyle{definition}
\newtheorem{rem}[thm]{Remark}

\theoremstyle{remark}

\begin{document}
\pagestyle{plain}
\today

\title{Metrics of constant curvature on a Riemann  surface with two corners
on the boundary}
\author{J\"urgen Jost}
\address{Max Planck Institute for Mathematics in the Sciences, Inselstr. 22, D-04013, Germany}
\email{jost@mis.mpg.de}

\author{Guofang Wang}
\address{University of Magdeburg,
Faculty of Mathematics,
P.O.Box 4120, D-39016 Magdeburg, Germany}
\email{guofang.wang@math.uni-magdeburg.de}

\author{Chunqin Zhou}
\address{Department of Mathematics, Shanghai Jiaotong University, Shanghai, 200240, China }
\email{cqzhou@sjtu.edu.cn}
\thanks{The third named author
supported partially by NSFC of China (No. 10301020)}

\maketitle

\section{Introduction}
On a Riemann surface, one of the  interesting geometric problems  is to determine which
functions can be realized as the Gaussian curvature of some
pointwise conformal metric. The classical uniformization theorem
tell us  that every smooth Riemannian metric on a two-dimensional
surface is pointwise conformal to one with constant curvature. This
question is by now well understood from many different perspectives,
and successfully approached by many different methods.

On this basis, research can  move on to  surfaces with
singularities. This, however, is by no means a straightforward
generalization of the smooth case. Results for smooth surfaces might not be true for
surfaces with singularities. For instance, there exist many
surfaces with conical singularities that do not admit a conformal
metric of constant Gauss curvature. In fact, a closed surface with
two conical singularities admits a conformal metric of constant
Gauss curvature if and only if its singularities have the same
angle and are in antipodal positions -- thus, such a surface
necessarily has the shape of an American football; this was proved by Troyanov \cite{T1}. Therefore a surface with exactly one singularity
(the teardrop) does not carry a conformal metric of constant Gauss
curvature.

This result was obtained by methods from complex analysis. It is
known, however, that the existence question for conformal metrics is
intimately linked to the Liouville equation. In recent years, very
powerful PDE methods have been developed to precisely determine the
asymptotic behavior of solutions of this equation near
singularities.

The purpose of the present paper then is to bring to
bear the full force of those methods on the existence problem for
conformal metrics with prescribed singularities. In fact, we shall
investigate the more general situation of surfaces with boundary. When
we have a boundary, the natural curvature condition there, the
analogue of the constant Gauss curvature condition in the interior, is
the one of constant geodesic curvature.

To continue the discussion about surfaces with singularities, let
us first recall  their definition, following  \cite{T1}. A conformal metric $ds^2$ on a Riemannian surface
$\Sigma$ (possibly with boundary) has a conical singularity of order
$\alpha$ (a real number with $\alpha>-1$) at a point $p\in
\Sigma\cup
\partial \Sigma$ if in some neighborhood of $p$
$$ds^2=e^{2u}|z-z(p)|^{2\alpha}|dz|^2$$
where $z$ is a coordinate of $\Sigma$ defined in this neighborhood
and $u$ is smooth away from $p$ and  continuous at $p$. The point
$p$ is then said to be a conical singularity of {\it angle}
$\theta=2\pi(\alpha+1)$ if $p\notin\partial \Sigma$ and a {\it
corner} of angle $\theta =\pi(\alpha +1)$ if $p\in \partial
\Sigma$. For example, a football has two singularities of equal
angle, while a teardrop has only one singularity. Both these
examples correspond to the case $-1<\alpha <0$; in case $\alpha >0$,
the angle is larger than $2\pi$, leading to a different geometric
picture. Such singularities also appear in orbifolds and branched
coverings. They can also describe the ends of complete Riemannian
surfaces with finite total curvature. If $(\Sigma, ds^2)$ has
conical singularities of order $\alpha_1, \alpha_2, \cdots,
\alpha_n$ at $p_1, p_2, \cdots, p_n$, then $ds^2$ is said to
represent the divisor {\bf A}$:= \Sigma^n_{i=1}\alpha_ip_i$.

For  a closed surface with more than two conical singularities, the existence problem of constant Gauss curvature
already becomes subtle. When all singularities have order
$\alpha\in (-1,0)$, Luo and Tian \cite{LT} gave a necessary and
sufficient condition. For the  case of general $\alpha$, a necessary and
sufficient condition was given by \cite{UY} recently for a closed
surface with 3 conical singularities. See also \cite{E} for a
simpler proof.

As already mentioned, the objective of this paper is to consider  surfaces (with
boundary) with corners on their boundary and to study the existence
problem of conformal metrics with constant Gauss curvature and
constant geodesic curvature on their boundary. Our first result
shows  that a disk with two corners admits a conformal metric
with constant Gauss curvature and constant geodesic curvature on
its boundary if and only if the two corners have the same angle.
This is analogous to the result of \cite{T2}. The disk is conformally
equivalent to ${\mathbb R}^2_+\cup\{\infty\}$. Note that the case of  a metric with zero geodesic curvature on
its boundary can be reduced to Troyanov's result.

\begin{thm} \label{thm1} It is possible to construct a  metric $g$
  with constant Gauss
curvature on the unit disk $D$ and constant geodesic curvature on
$\Gamma_\pm:=\partial D\cap \{(x,y)\in \R^2\,|\,\pm y>0\}$
admitting two corners  $p_1=(1,0)$ with order $\alpha_1>-1$ and
$p_2=(-1,0)$ with order $\alpha_2>-1$  if and
only if
\[\alpha_1=\alpha_2.\]
\end{thm}

In Theorem \ref{thm1}, the constant geodesic curvatures on
$\Gamma_+$ and $\Gamma_-$  may be different. All solutions can
be explicitly written down, see Theorem \ref{theo-metr}. Theorem
\ref{thm1} is not difficult to prove. But it is a good starting point
for our research.

What we do in fact is more general than this generalization of Troyanov's result.
 Let us  denote $\R^{2}_{+}=\{(s,t)|t>0\}$.  We consider

\beq \label{eq-1}
\left\{
\begin{array}{rcll}
-{\Delta} u &=& \ds\vs |x|^{2\alpha}e^u, &\qquad\text{in }
 \R^{2}_{+},\\
\ds\vs \frac {\partial u}{\partial t}&=& \ds c_1 |x|^{\alpha}e^{\frac u2},
& \qquad \text{on } \partial \R^{2}_{+}\cap\{s>0\}\\
\ds\vs\frac {\partial u}{\partial t}&=&\ds c_2 |x|^{\alpha}e^{\frac u2}, &
\qquad \text{on }\partial \R^{2}_{+}\cap\{s<0\}\\
\end{array}
\right. \eeq with the energy conditions

\beq \label{con-1} \int_{\R^{2}_{+}}|x|^{2\alpha}e^udx<\infty,
 \int_{\partial \R^{2}_{+}}|x|^{\alpha}e^{\frac
u2}ds<\infty. \eeq
Here $c_1,$ $c_2$ are  constants and $\alpha>-1$.

We call $u\in
H^{1}_{loc}(\overline{\R^2_+})$  a weak solution of
(\ref{eq-1})-(\ref{con-1}) if it satisfies
$$
\int_{\R^2_+}\nabla u\cdot\nabla\varphi
dx+c_1\int_{\partial\R^2_+\cap\{s>0\}}|x|^{\alpha}e^{\frac
u2}\varphi
ds+c_2\int_{\partial\R^2_+\cap\{s<0\}}|x|^{\alpha}e^{\frac
u2}\varphi ds=\int_{\R^2_+}|x|^{2\alpha}e^u\varphi dx
$$
for any smooth function $\varphi(x)$ on $\overline{\R^2_+}$ with compact
support.
% Here
%$H^{1}_{loc}(\overline{\R^2_+}) = \{u\in H^1(\Omega):  \Omega
%\subset \overline{\R^2_+}\}$ and  $\Omega$ is any bounded subset.
Since $u\in H^{1}_{loc}(\overline{\R^2_+})$ implies $e^u\in
L^p_{loc}(\overline{\R^2_+})$ for all $p>1$,  by standard
elliptic regularity we conclude that any weak solution $u$ of
(\ref{eq-1}) is a classical solution when $\alpha\geq 0$ while   $u$
is smooth away from the origin and $u\in W^{2,q}$ near the origin
for $1<q<-\frac 1\alpha$ when $-1<\alpha <0$. In particular, $u$
is continuous at the origin in any case. In the sequel, we assume
that a solution $u$ of (\ref{eq-1})-(\ref{con-1})
always satisfies $u\in C^2(\R^{2}_{+})\cap
C^1(\overline{\R^{2}_{+}}\setminus\{0\})$ and that $u$ is continuous at
the origin.

Geometrically,  a solution $u$ of (\ref{eq-1})
-(\ref{con-1}) determines a metric
$ds^2=|z|^{2\alpha}e^u|dz|^2$ with constant scalar curvature
$1$ on $\R^{2}_{+}$ and with geodesic curvature $-c_1$ on $\partial
\R^{2}_{+}\cap\{s>0\}$ and geodesic curvature $-c_2$ on $\partial
\R^{2}_{+}\cap\{s<0\}$. Moreover $ds^2=|z|^{2\alpha}e^u|dz|^2$ has
a conical singularity at $z=0$. Let  $1$
and $-1$ be two points on the boundary of  the unit disk $D$ . We take a
conformal transformation $\phi$ mapping $D$ to $\R^{2}_{+}$ and
$\partial D$ to $\partial \R^{2}_{+}$ with $\phi (1)=0$ and
$\phi(-1)=\infty$. With such a conformal transformation, the
metrics studied in Theorem \ref{thm1} are solutions of
(\ref{eq-1})-(\ref{con-1}). Our main result in this paper is to
show the converse, namely, any solution of
(\ref{eq-1})-(\ref{con-1})  is in fact obtained from a metric in
Theorem \ref{thm1}.

\begin{thm}\label{theo-metr} Let $u$ be a solution of (\ref{eq-1})-(\ref{con-1}).
Then
 $ds^2=e^{u}|z|^{2\alpha} |dz|^2$ comes from a conformal
metric as in Theorem \ref{thm1}.  More precisely, there exists
$\lambda
>0$ such that:
\begin{enumerate}
\item[(1)]When $\alpha=2k, k=0,1,2,\dots$, then $c_1=c_2$.
And when $\alpha=2k+1, k=0,1,2,\dots$, then $c_1=-c_2$. In this case the metric is
$$
ds^2=\frac{8(\alpha+1)^2\lambda^{2\alpha+2}|z|^{2\alpha}|dz|^2}{(\lambda^{2\alpha
+2}+|z^{\alpha+1}-z_0|^2)^2}$$ for some $z_0=(s_0,t_0)$ with
$s_0\in \R$ and $t_0=\frac{c_1\lambda^{\alpha+1}}{\sqrt{2}}$.
\item[(2)] When $\alpha\neq k, k=0,1,2,\dots$, then for any $c_1$ and $c_2$, the metric is
$$
ds^2=\frac{8(\alpha+1)^2\lambda^{2\alpha+2}|z|^{2\alpha}|dz|^2}{(\lambda^{2\alpha
+2}+|z^{\alpha+1}-z_0|^2)^2}$$ for some $z_0=(s_0,t_0)$ with
$s_0=\frac{\lambda^{\alpha+1}(c_1\cos(\pi\alpha)-c_2)}{\sqrt{2}\sin(\pi\alpha)}$
and $t_0=\frac{c_1\lambda^{\alpha+1}}{\sqrt{2}}$.
\end{enumerate}
\end{thm}

This result  is a natural generalization of the classification
result of Chen-Li \cite{CL2} for the Liouville equation
\begin{equation}\label{L}
-{\Delta} u =e^u \quad \hbox { in } \R^2\end{equation} with finite
area $\int_{\R^2}e^u<\infty$ and  the classification result of
Li-Zhu \cite{LZ} for solutions of
\begin{equation}\label{L2}
\left \{
\begin{array}{rcll}
-{\Delta} u & =& \ds\vs e^u & \quad \hbox { in } \R^2_+ \\
 \ds\frac {\partial u}{\partial t}&=& \ds c e^{\frac u 2} &\quad\hbox {on }\partial \R^2_+.
 \end{array}
 \right.
\end{equation}

Geometrically,  the result of Chen-Li covers the case of the
standard sphere. In fact, their classification result tells us
that any solution of the Liouville equation (\ref{L}) with finite
area can be compactified as a metric on the standard sphere with
constant curvature. Similarly, the result of Li-Zhu deals with a
portion of the standard sphere cut by a 2-plane. Namely, from
their result we know that any solution of (\ref{L2}) can be
compactified as a metric on such a portion of the standard sphere
with constant Gauss curvature and constant geodesic curvature on
the boundary. In this spirit, our result (for $-1<\alpha<1$) then
deals with a portion of the standard sphere cut by two 2-planes
with angle $\pi(\alpha+1)$.

%\setlength{\unitlength}{1mm}
%\begin{center}
%\scalebox{0.6}[0.6]{\includegraphics[0,0][180, 180]{picture.eps}}
%\end{center}

It would then be interesting to consider portions of the
standard sphere cut by 3 or more 2-planes. This is related to the
result of Umehara-Yamada \cite{UY}, see also \cite{E}. We will
return to this issue later. In another direction, our result is  a
generalization of Prajapat-Tarantello \cite{PT}, who classify
solutions of the Liouville equation with one singularity. For the case
$c_1=c_2=0$,  Theorem \ref{theo-metr}  can
be reduced to their result. For  other classication results, or
different proofs, see \cite{CL3}, \cite{CW},\cite{CY},\cite{HT},
\cite{HW}, \cite{JW} , \cite{M},  and \cite{Z}.

Our method to deal with (\ref{eq-1})-(\ref{con-1}) can be viewed
as a combination of the methods developed for those previous
results. We shall make particular use of \cite{JLW} and \cite{LZ}.
The main issue is the determination of
$$
d:=-\lim_{|x|\rightarrow \infty}\frac{u(x)}{\ln|x|}.
$$
Note that equations (\ref{eq-1}) are no longer translation
invariant and  a
solution of (\ref{eq-1})-(\ref{con-1}) will no longer be radially symmetric  if one of $c_i\neq 0$ for
$i=1,2$. The methods used in \cite{LZ} and \cite{PT} can therefore not
be directly utilized to prove Theorem \ref{theo-metr}. However, after we have
shown that the metric
$ds^2=e^{\tilde{u}}|dz|^2=|z|^{2\alpha}e^u|dz|^2$ has two conical
singularities at $z=0$ and $z=\infty$, we can define
$$
\eta(z)=\left(\frac{\partial^2 \tilde{u}}{\partial z^2}-\frac
12(\frac{\partial\tilde{u}}{\partial z})^2\right)|dz|^2,\quad z \text {
in } \R^2_+,
$$ which can be  extended to a projective connection on $\S^2=\C\cup
\{\infty\}$
as defined in \cite{T2}. Then the problem is reduced to a linear partial
differential system, see (\ref{eq-10}) and (\ref{con-5}).
Finally we solve this boundary problem and demonstrate Theorem
\ref{theo-metr}.

\
\

\section{Projective Connections}
In this section, we will state the definition and the properties of
the  projective connection discussed in the papers
of Troyanov \cite{T2} and Mandelbaum \cite{Ma}. In the last
section, we will demonstrate our main result in the sense of a
projective connection on $\C\cup \{\infty\}$.

\
\

Assume that $\Sigma$ is a Riemann surface. Let $\eta $ be a
quadratic differential. If

\quad   (1) $\eta(z)=\phi(z)|dz|^2$ is a meromorphic quadratic
differential in each local coordinate $(U, z)$ on $\Sigma$,

 \quad
(2)$ \eta(w)=\eta(z)+\{z,w\}|dw|^2$ in the overlap of two local
coordinates $(U,z)$ , $(V,w)$,\\
then $\eta$ is called  a {\bf projective connection} on $\Sigma$.
Here $\{,\}$ denotes the Schwarzian derivative:
$$
\{z,w\}=\frac{z'''}{z'}-\frac 32(\frac{z''}{z'})^2
$$  a function $z$ of $w$.

\
\

A point $p\in\Sigma$ is called a regular point of the projective
connection $\eta$ if $\eta $ is holomorphic at this point. We say
that $\eta$ has a {\bf regular singularity} of weight $\rho$ at
$p$ if $$ \eta(z)=(\frac {\rho}{z^2}+\frac
{\sigma}{z}+\phi_1(z))|dz|^2$$ where $\phi_1(z) $ is holomorphic,
and $z$ is a local coordinate at $p$ with $z(p)=0$.

\
\

The projective connection is said to be compatible with the
divisor {\bf A }:$=\Sigma^{n}_{i=1}\alpha_ip_i$ if it is
regular in $\Sigma-\{p_1,\cdots , p_n\}$ and has, for each $i$, a
regular singularity of weight $\rho_i=-\frac
12\alpha_i(\alpha_i+2)$ at $p_i$. The next two lemmas are examples
about some results for the projective connection from \cite{T2}.

\begin{lm}
The definition of the weight for a singular point is independent
of the choice of local coordinate.
\end{lm}

\begin{lm}
If $ds^2=e^{u}|dz|^2$ is a conformal metric of constant curvature on
$\Sigma$ representing the divisor {\bf A} then $$
\eta(z)=(\frac{\partial^2u}{\partial z^2}-\frac 12(\frac{\partial
u}{\partial z})^2)|dz|^2$$ defines a projective connection
compatible with the divisor {\bf A}.
\end{lm}

\
\

\section{Asymptotic behavior}
We will first rewrite the equation (\ref{eq-1}). Set
$\widetilde{u}=u+2\alpha \ln{|x|}$. Then $\widetilde{u}$ satisfies

\beq \label{eq-3} \left\{
\begin{array}{rcll}
-{\Delta} \widetilde{u} &=& e^{\widetilde{u}}, &\qquad \text{in}
\quad \R^{2}_{+},\\
\frac {\partial \widetilde{u}}{\partial t}&=& c_1 e^{\frac
{\widetilde{u}}{2}}, &\qquad \text{on}\quad
\partial \R^{2}_{+}\cap\{s>0\},\\
\frac {\partial \widetilde{u}}{\partial t}&=& c_2 e^{\frac
{\widetilde{u}}{2}}, &\qquad \text{on}\quad
\partial \R^{2}_{+}\cap\{s<0\}
\end{array}
\right. \eeq with the energy conditions

\beq \label{con-3} \int_{\R^{2}_{+}}e^{\widetilde{u}}dx<\infty ,
\eeq

\beq \label{con-4} \int_{\partial \R^{2}_{+}}e^{\frac
{\widetilde{u}}{2}}ds<\infty. \eeq

\begin{prop}\label{prop-bdth}
Any solution $\widetilde{u}$   of (\ref{eq-3}) with
(\ref{con-3}) and (\ref{con-4}) is bounded from above in the region $\overline{\R^{2}_{+}\setminus
B^{+}_{\varepsilon}(0)}$, for each $\varepsilon
>0$.
\end{prop}

To prove Proposition \ref{prop-bdth}, we need the following Lemma.
\begin{lm}\label{lm-11}
Assume that $u$ is a solution of

\begin{equation*}
\left\{
\begin{array}{rcll}
-{\Delta} u &=& 0, &\qquad \text{in}
\quad B^{+}_{R},\\
\frac {\partial u}{\partial t}&=& f(x), &\qquad \text{on}\quad
\{t=0\}\cap \partial B^{+}_{R}\\
u&=& 0, &\qquad \text{on}\quad \partial B^{+}_{R}\cap
\overline{B}^{+}_{R}.
\end{array}
\right.
\end{equation*}

\noindent with $f\in L^1(\{t=0\}\cap \partial B^{+}_{R})$ for any
$R>0$. Then for every $\delta_1\in (0,4\pi)$ we have
$$
\int_{B^{+}_{R}}exp\{\frac
{(4\pi-\delta_1)|u(x)|}{||f||_1}\}dx\leq\frac{16\pi^2R^2}{\delta_1}
$$
and for every $\delta_2\in (0,2\pi)$
$$
\int_{\partial B^{+}_{R}\cap \{t=0\}}exp\{\frac
{(2\pi-\delta_2)|u(x)|}{||f||_1}\}ds\leq\frac{4\pi R}{\delta_2}
$$
where $||f||_1=\int_{\{t=0\}\cap \partial B^{+}_{R}}|f|ds$.

\end{lm}

\begin{proof}
Set

$$\Gamma_1=\{t=0\}\cap \overline{B}^{+}_{R}, \quad
\Gamma_2=\{t>0\}\cap \partial B^{+}_{R}.$$

\noindent Let
$$
\phi(y)=\frac
{1}{2\pi}\int_{\Gamma_1}(\log{\frac{2R}{|x-y|}}+\log{\frac{2R}{|x-\overline{y}|}})|f(x)|dx
$$
where $\overline {y}$ is the reflection point of $y$ about
$\{t=0\}$.

\noindent A direct computation yields

\begin{equation*}  \left\{
\begin{array}{rcll}
-{\Delta} \phi &=& 0, &\qquad \text{in}
\quad B^{+}_{R},\\
\frac {\partial \phi}{\partial t}&=& -|f|,
&\qquad \text{on}\quad \Gamma_1.\\
\end{array}
\right. \end{equation*}

\noindent Note that $\phi \geq 0$ for $x\in B^{+}_{R}$ since
$\frac{2R}{|x-y|}\geq 1$ for any $x,y\in B^{+}_{R}$. We have

\begin{equation*}  \left\{
\begin{array}{rcll}
-{\Delta} (u-\phi) &=& 0, &\qquad \text{in}
\quad B^{+}_{R},\\
\frac {\partial(u-\phi)}{\partial t}&=& f+|f|,
&\qquad \text{on}\quad \Gamma_1\\
u-\phi &\leq & 0,&\qquad \text{on} \quad\Gamma_2.
\end{array}
\right. \end{equation*} It follows from the maximum principle and the
Hopf Lemma that $u\leq \phi$ in $\overline{B}^{+}_{R}$.

\noindent By a similar argument we also have

\begin{equation*}  \left\{
\begin{array}{rcll}
-{\Delta} (u+\phi) &=& 0, &\qquad \text{in}
\quad B^{+}_{R},\\
\frac {\partial(u+\phi)}{\partial t}&=& f-|f|,
&\qquad \text{on}\quad \Gamma_1\\
u+\phi &\geq & 0,&\qquad \text{on} \quad\Gamma_2.
\end{array}
\right. \end{equation*} which implies that $u\geq -\phi$ in
$\overline{B}^{+}_{R}$. Therefore we have $|u|\leq \phi$ in
$\overline{B}^{+}_{R}$ and thus we have
$$
\int_{B^{+}_{R}}exp\{\frac
{(4\pi-\delta_1)|u(x)|}{||f||_1}\}dx\leq
\int_{B^{+}_{R}}exp\{\frac {(4\pi-\delta_1)\phi}{||f||_1}\}dx.
$$
and
$$
\int_{\Gamma_1}exp\{\frac {(2\pi-\delta_2)|u(x)|}{||f||_1}\}ds\leq
\int_{\Gamma_1}exp\{\frac {(2\pi-\delta_2)\phi}{||f||_1}\}ds.
$$

At this point, using  Jensen's inequality, we can follow  the argument
of \cite{BM}, proof of Theorem 1, to
conclude the result.

\end{proof}

{\bf Proof of Proposition \ref{prop-bdth}} We first fix
$\varepsilon>0$, and assume that $\widetilde{u}$ is  a solution
of (\ref{eq-3}) with (\ref{con-3}) and (\ref{con-4}). From Theorem
2 of \cite{BM} it suffices to show that, for any $x_0 \in \partial
\R^{2}_{+}\setminus \overline{B}^{+}_{\varepsilon}(0)$,
$\widetilde{u} $ is bounded from above on
$\overline{B}^{+}_{R}(x_0)$ for some small number $R>0$, with a
bound that is independent of the point $x_0$. In the following, we
denote by $C$ various constants independent of $x_0$.

Write $g=e^{\widetilde{u}}$, $f=c(x)e^{\frac {\widetilde{u}}{2}}$
where $c(x)$ is a function on $\partial \R^2_+\setminus \{0\}$
with $c(x)=c_1$ when $s>0$ and $c(x)=c_2$ when $s<0$, where we
write $x=(s,t)$. Then $\widetilde{u}$ satisfies
\begin{equation*}  \left\{
\begin{array}{rcl}
-{\Delta} \widetilde{u} &=& g, \qquad \text{in}
\quad B^+_{R}(x_0),\\
\frac{\partial \widetilde{u}}{\partial t}&=& f,
\qquad\text{on}\quad \Gamma_1.
\end{array}
\right. \end{equation*}

It is clear that $f\in L^1(\partial \R^2_+)$. Set $f =f_1+f_2$
with $||f_1||_{L^1(\partial \R^{2}_{+})}\leq \pi$ and $f_2\in
L^{\infty}(\partial \R^2_+)$. Let $\Gamma_1$ and $\Gamma_2$ be as
Lemma \ref{lm-11}. Define $\widetilde{u}_1$, $\widetilde{u}_2$ and
$\widetilde{u}_3$ by

\begin{equation*}  \left\{
\begin{array}{rcll}
-{\Delta} \widetilde{u}_1 &=& e^{\widetilde{u}}, &\qquad
\text{in}
\quad B^{+}_{R}(x_0),\\
\frac {\partial \widetilde{u}_1}{\partial t}&=& 0,
&\qquad \text{on}\quad \Gamma_1\\
\widetilde{u}_1&=& 0, &\qquad\text{on}\quad \Gamma_2.
\end{array}
\right. \end{equation*}

\begin{equation*}  \left\{
\begin{array}{rcll}
-{\Delta} \widetilde{u}_2 &=& 0, &\qquad \text{in}
\quad B^{+}_{R}(x_0),\\
\frac {\partial \widetilde{u}_2}{\partial t}&=& f_1,
&\qquad \text{on}\quad \Gamma_1\\
\widetilde{u}_2&=& 0, &\qquad\text{on}\quad \Gamma_2.
\end{array}
\right. \end{equation*}

\begin{equation*}  \left\{
\begin{array}{rcll}
-{\Delta} \widetilde{u}_3 &=& 0, &\qquad \text{in}
\quad B^{+}_{R}(x_0),\\
\frac {\partial \widetilde{u}_3}{\partial t}&=& f_2,
&\qquad \text{on}\quad \Gamma_1\\
\widetilde{u}_3&=& 0, &\qquad\text{on}\quad \Gamma_2.
\end{array}
\right. \end{equation*}

\noindent Extending $\widetilde{u}_1$ evenly we have
\begin{equation*}  \left\{
\begin{array}{rcll}
-{\Delta} \widetilde{u}_1 &=& e^{\widetilde{u}}, &\qquad
\text{in}
\quad B_{R}(x_0),\\
\widetilde{u}_1&=& 0, &\qquad\text{on}\quad \partial B_R(x_0).
\end{array}
\right. \end{equation*} By using Theorem 2 in \cite{BM} and
(\ref{con-3}) we have
$$
||\widetilde{u}_1||_{L^{\infty}(\overline{B}^+_R(x_0))}\leq C.
$$
For $\widetilde{u}_2$, by Lemma \ref{lm-11}, we have
$$
\int_{B^+_R(x_0)}\exp(2|\widetilde{u}_2|)dx\leq C,\qquad
\int_{\Gamma_1}\exp(|\widetilde{u}_2|)ds\leq C
$$
and in particular $||\widetilde{u}_2||_{L^q(B^+_R(x_0))}\leq C$
and $||\widetilde{u}_2||_{L^q(\Gamma_1)}\leq C$ for any $q>1$.

\noindent For $\widetilde{u}_3$, it is obvious that
$$
||\widetilde{u}_3||_{L^{\infty}(\overline{B}^+_{\frac
R2}(x_0))}\leq C.
$$

\noindent Let
$\widetilde{u}_4=\widetilde{u}-\widetilde{u}_1-\widetilde{u}_2-\widetilde{u}_3$.
Then we have
\begin{equation*}  \left\{
\begin{array}{rcl}
-{\Delta} \widetilde{u}_4 &=& 0, \qquad \text{in}
\quad B^+_{R}(x_0),\\
\frac{\partial\widetilde{u}_4}{\partial t}&=& 0,
\qquad\text{on}\quad \Gamma_1.
\end{array}
\right. \end{equation*} Extending $\widetilde{u}_4$ evenly,
$\widetilde{u}_4$ becomes a harmonic function on $B_R(x_0)$. Then
the mean
value theorem for harmonic functions implies that
$$
||\widetilde{u}^+_4||_{L^{\infty}(\overline{B}^+_{\frac
R2}(x_0))}\leq C ||\widetilde{u}^+_4||_{L^{1}(B^+_R(x_0))}.
$$
Notice that
$$
\widetilde{u}^+_4\leq
\widetilde{u}^++|\widetilde{u}_1|+|\widetilde{u}_2|+|\widetilde{u}_3|,
$$
and
$$
\int_{\R^2_+}\widetilde{u}^+dx\leq
\int_{\R^2_+}e^{\widetilde{u}^+}dx<\infty.
$$
We get
$$
||\widetilde{u}^+_4||_{L^{\infty}(\overline{B}^+_{\frac
R2}(x_0))}\leq C.
$$

\noindent Finally, we write
\begin{equation*}  \left\{
\begin{array}{rcl}
-{\Delta} \widetilde{u} &=& e^{\widetilde{u}}=g, \qquad \text{in}
\quad B^+_{R}(x_0),\\
\frac{\partial \widetilde{u}}{\partial t}&=&
c(x)e^{\frac{\widetilde{u}}{2}}=f, \qquad\text{on}\quad \Gamma_1.
\end{array}
\right. \end{equation*} The standard elliptic estimates imply that
$$
||\widetilde{u}^+||_{L^{\infty}(\overline{B}^+_{\frac
R4}(x_0))}\leq C ,
$$
since $||f||_{L^q(B^+_{\frac R2}(x_0))}\leq C$ and
$||g||_{L^q(\partial B^+_{\frac R2}(x_0))\cap \{t=0\})}\leq C$ for
any $q>1$.\qed

\

\

As in the proof of Proposition \ref{prop-bdth}, in the sequel we
always let $c(x)$ be a function on $\partial \R^2_+\setminus\{0\}$
with $c(x)=c_1$ when $s>0$ and $c(x)=c_2$ when $s<0$, where
$x=(s,t)$. In virtue of Proposition \ref{prop-bdth}, we obtain the
asymptotic behavior of the solution of (\ref{eq-1})-(\ref{con-1}).
More precisely, we have the following

\begin{prop}\label{prop-ayth}
Let $u$ be a solution of (\ref{eq-1})-(\ref{con-1}). Let
$$
d=\frac{1}{\pi} \int_{\R^{2}_{+}}|x|^{2\alpha}e^udx-\frac{1}{\pi}
\int_{\partial \R^{2}_{+}}c(x)|x|^{\alpha}e^{\frac u2}ds.
$$
Then we have
$$
\lim_{|x|\rightarrow \infty}\frac{u(x)}{\ln|x|}=-d.
$$
\end{prop}

\begin{proof}
Let \begin{eqnarray*} w(x) &=&\frac{1}{2\pi}
\int_{\R^{2}_{+}}(\log |x-y|+\log
|\overline{x}-y|-2\log|y|)|y|^{2\alpha}e^{u(y)}dy\\
& &-\frac{1}{2\pi} \int_{\partial \R^{2}_{+}}(\log |x-y|+\log
|\overline{x}-y|-2\log|y|)c(y)|y|^{\alpha}e^{\frac {u(y)}{2}}dy.
\end{eqnarray*}
where $\bar{x}$ is the reflection point of $y$ about $\{t=0\}$. It
is easy to check that $w(x)$ satisfies

\begin{equation*}  \left\{
\begin{array}{rcll}
{\Delta} w &=& |x|^{2\alpha}e^u, &\qquad \text{in}
\quad \R^{2}_{+},\\
\frac {\partial w}{\partial t}&=& -c(x) |x|^{\alpha}e^{\frac u2},
&\qquad \text{on}\quad \partial \R^{2}_{+}\setminus\{0\}.
\end{array}
\right. \end{equation*}
and
$$
\lim_{|x|\rightarrow \infty}\frac{w(x)}{\ln|x|}=d.
$$
Consider $v(x)=u+w$. Then $v(x)$ satisfies
\begin{equation*}\left\{
\begin{array}{rcl}
{\Delta} v &=& 0, \qquad \text{in}
\quad \R^{2}_{+},\\
\frac {\partial v}{\partial t}&=& 0, \qquad \text{on}\quad
\partial \R^{2}_{+} \setminus\{0\}.
\end{array}
\right. \end{equation*}

We extend $v(x)$ to $\R^2$ by even reflection such that $v(x)$ is
harmonic in $\R^2$ . From Proposition \ref{prop-bdth} we know
$v(x)\leq C(1+\ln (|x|+1))$ for some positive constant $C$. Thus
$v(x)$ is a constant. This completes the proof.
\end{proof}

\begin{rem}\label{rem-1}From (\ref{con-1}), it is  easy to check that $d\geq 2+2\alpha$.
\end{rem}

\
\

\section{The exact value of $d$}
In this section, we want to compute the value of $d$. We need to
distinguish
two cases. When $c_1\leq 0$ and $c_2\leq 0$, we will employ a
similar argument as in \cite{JLW} when they proved $\gamma_i <2$ in
proposition 7.1 to show that $d>2+2\alpha$. Here $c_1\leq 0$ and
$c_2\leq 0$ are crucial such that $w(x)<0$ in $D^+$, see
Proposition \ref{prop-dth1}. Once we have proved that
$d>2+2\alpha$, we can obtain an extension of $u(x)$ near $\infty$,
see (\ref{asy-est4}). Then we can use the Pohozaev identity of
(\ref{eq-1}) to prove $d=4+4\alpha$. When $c_i>0$ for $i=1$ or
$i=2$, this method will not work well. We will use the moving
sphere method , which was used in \cite{LZ}, to show
$d>2(1+\alpha)$. Let us start with the negative case.

\begin{prop}\label{prop-dth1}If  $c_1\leq 0$ and
$c_2\leq 0$ in (\ref{eq-1})-(\ref{con-1}), we have $d>2+2\alpha$.
\end{prop}

\begin{proof}
Assume by contradiction that $d=2+2\alpha$. Let $v$ be the Kelvin
transformation of $u$, i.e. $v(x)=u(\frac
{x}{|x|^2})-(4\alpha+4)\ln|x|$. Then $v$ satisfies

\begin{equation*}  \left\{
\begin{array}{rcl}
-{\Delta} v &=& |x|^{2\alpha}e^v, \qquad \text{in}
\quad \R^{2}_{+},\\
\frac {\partial v}{\partial t}&=& c(x) |x|^{\alpha}e^{\frac v2},
\qquad \text{on}\quad \partial \R^{2}_{+}\setminus\{0\}.
\end{array}
\right. \end{equation*} with the energy conditions
$$
\int_{\R^{2}_{+}}|x|^{2\alpha}e^vdx<\infty. $$ and $$
\int_{\partial \R^{2}_{+}}|x|^{\alpha}e^{\frac v2}dt<\infty.
$$
Here $c(x)$ is a function as in the above section.

Let $D^+$ be a small half disk centered at zero. Define $w(x)$ by
\begin{eqnarray*} w(x) &=&\frac{1}{2\pi}
\int_{D^+}(\log |x-y|+\log
|\overline{x}-y|)|y|^{2\alpha}e^{v(y)}dy\\
& &-\frac{1}{2\pi} \int_{\partial D^+\cap\{t=0\}}(\log |x-y|+\log
|\overline{x}-y|)c(y)|y|^{\alpha}e^{\frac {v(y)}{2}}dy.
\end{eqnarray*}
and define  $g(x)=v(x)+w(x)$. It is clear that
\begin{equation*}\left\{
\begin{array}{rcl}
{\Delta} g &=& 0, \qquad \text{in}
\quad D^+,\\
\frac {\partial g}{\partial t}&=& 0, \qquad \text{on}\quad
\{\partial D^+\cap\{t=0\}\}\setminus\{0\} .
\end{array}
\right. \end{equation*} Therefore by extending $g(x)$ to
$D\setminus\{0\}$ evenly we obtain a harmonic $g(x)$ in
$D\setminus\{0\}$.

On the other hand, we can check that
$$
\lim_{|x|\rightarrow 0}\frac {w}{-\log{|x|}}=0
$$
which implies
$$
\lim_{|x|\rightarrow 0}\frac
{g(x)}{-\log{|x|}}=\lim_{|x|\rightarrow 0}\frac
{v(x)+w(x)}{-\log{|x|}}=2\alpha +2.
$$
Since $g(x)$ is harmonic in $D\backslash\{0\}$, we have
$g(x)=-(2\alpha +2)\log{|x|}+g_0(x)$ with a smooth harmonic
function $g_0$ in $D$. By the definition, we have $w(x)<0$ since
$c(x)$ is negative. Thus, we have
$$
\int_{D^+}|x|^{2\alpha}e^{v}dx=\int_{D^+}|x|^{2\alpha}e^{g-w}dx\geq
\int_{D^+}|x|^{2\alpha}|x|^{-2\alpha-2}e^{g_0}dx=\infty,
$$
which is a contradiction with
$\int_{\R^2_+}|x|^{2\alpha}e^{v}dx<\infty$. Hence we have shown
that $d> 2\alpha+2$.

\end{proof}

From $d> 2\alpha+2$ we can improve the estimates for $e^u$ to
\begin{equation}\label{asy-est1}
e^{u}\leq C|x|^{-2-2\alpha-\varepsilon}, \qquad \text{for}\quad
|x| \quad \text{near} \quad \infty.
\end{equation}
Then by using  potential analysis,  we  obtain
$$
-d\ln{|x|}-C\leq u(x)\leq -d\ln{|x|}+C
$$
for some constant $C>0$ and $\varepsilon>0$, see \cite{CL2}.
Furthermore following the  idea for the
derivation of gradient estimates in \cite{CK} and \cite{WZ}, we
get
\begin{equation*}
|\langle x, \nabla u \rangle+d|\leq C|x|^{-\varepsilon} \qquad
\text{for} \quad |x|\quad \text{near}\quad \infty,
\end{equation*}
consequently we have
\begin{equation}\label{asy-est2}
|u_r+\frac{d}{r}|\leq C|x|^{-1-\varepsilon} \qquad \text{for}
\quad |x|\quad \text{near}\quad \infty.
\end{equation}
In a similar way, we can also get
\begin{equation}\label{asy-est3}
|u_{\theta}|\leq C|x|^{-\varepsilon} \qquad \text{for} \quad
|x|\quad \text{near}\quad \infty.
\end{equation}
From (\ref{asy-est2}) and (\ref{asy-est3}) we can also get by
 standard potential analysis that
\begin{equation}\label{asy-est4}
u(x)=-d\ln{|x|}+C+O(|x|^{-1}) \qquad \text{for}\quad |x| \quad
\text{near}\quad \infty,
\end{equation}
Here $(r,\theta)$ is the polar coordinate system on $\R^2$, and $C,
\varepsilon$ are some positive constants.

\begin{prop}\label{prop-dth2}If $d>2+2\alpha$, then we have
$d=4+4\alpha$.
\end{prop}

\begin{proof}
Firstly we establish the Pohozaev identity of
(\ref{eq-1})-(\ref{con-1}). Multiply equation (\ref{eq-1}) by
$x\cdot \nabla u$ and integrate over $B^+_R$ to obtain
$$
-\int_{B^+_R}(x\cdot \nabla u){\Delta}
udx=\int_{B^+_R}|x|^{2\alpha}e^ux\cdot \nabla udx
$$
Since
$$
(x\cdot \nabla u){\Delta} u=div((x\cdot \nabla u)\nabla
u)-div(\frac {x|\nabla u|^2}{2}),
$$
$$
|x|^{2\alpha}e^ux\cdot \nabla u=div(x|x|^{2\alpha
}e^u)-div(x)|x|^{2\alpha}e^u-e^ux\cdot \nabla |x|^{2\alpha},
$$
and
$$
x\cdot \nabla |x|^{2\alpha}=2\alpha |x|^{2\alpha},
$$
 we obtain
\begin{eqnarray*}
& &\int_{\partial B^+_R\cap \{t>0\}}x\cdot \nu \frac {|\nabla
u|^{2}}{2}-(\nu \cdot \nabla u)(x\cdot \nabla u)ds\\
& &+\int_{\partial B^+_R\cap \{t=0\}}x\cdot \nu \frac {|\nabla
u|^{2}}{2}-(\nu \cdot \nabla u)(x\cdot \nabla u)ds\\
& = & \int_{\partial B^+_R\cap \{t>0\}}x\cdot \nu
|x|^{2\alpha}e^uds+\int_{\partial B^+_R\cap \{t=0\}}x\cdot \nu
|x|^{2\alpha}e^uds \\
& &-(2+2\alpha)\int_{B^+_R}|x|^{2\alpha}e^udx,
\end{eqnarray*}
where $\nu$ is the outward unit normal vector to $\partial B^+_R$.
Hence we have
\begin{eqnarray*}
& & R\int_{\partial B^+_R\cap \{t>0\}}\frac {|\nabla u|^2}2-|\frac
{\partial u}{\partial r}|^2ds+\int_{\partial B^+_R\cap
\{t=0\}}\frac {\partial u}{\partial t}(x\cdot \nabla u)ds\\
&=& R\int_{\partial B^+_R\cap
\{t>0\}}|x|^{2\alpha}e^uds-(2+2\alpha)\int_{B^+_R}|x|^{2\alpha}e^udx.
\end{eqnarray*}
Since
\begin{eqnarray*}
& & \int_{\partial B^+_R\cap \{t=0\}}\frac {\partial u}{\partial
t}(x\cdot \nabla u)ds\\
&=& \int^R_{-R}c(x)|x|^{\alpha}e^{\frac u2}s{\partial}_s  uds\\
&=& 2\int^R_{-R}c(x)|x|^{\alpha}s\partial_s e^{\frac u2}ds\\
&=& 2c(x)|s|^\alpha se^{\frac
u2}|^R_{-R}-(2+2\alpha)\int^R_{-R}c(x)|x|^{\alpha}e^{\frac u2}ds,
\end{eqnarray*}
 we get the Pohozaev identity
\begin{eqnarray*}
& & R\int_{\partial B^+_R\cap \{t>0\}}\frac
{|u_{\theta}|^2}{2R^2}-\frac
{|u_r|^2}{2}ds\\
&=& R\int_{\partial B^+_R\cap
\{t>0\}}|x|^{2\alpha}e^uds-(2+2\alpha)\int_{B^+_R}|x|^{2\alpha}e^udx\\
& & -2c(x)|s|^\alpha se^{\frac
u2}|^R_{-R}+(2+2\alpha)\int^R_{-R}c(x)|x|^{\alpha}e^{\frac u2}ds.
\end{eqnarray*}
In virtue of (\ref{asy-est1}), (\ref{asy-est2}) and
(\ref{asy-est3}), we let $R\rightarrow \infty$ in the Pohozaev
identity and get
$$
d=4+4\alpha.
$$
\end{proof}

\
\

Next let us consider the case $c_i>0$ for $i=1$ or $i=2$.

\begin{prop}\label{prop-dth3}
If $c_i>0$ for $i=1$ or $i=2$, then $d\geq 4+4\alpha$ and
consequently $d=4+4\alpha$.
\end{prop}

\begin{proof}Without loss of  generality, we assume that
$c_1>0$. First we have $d\geq 2(1+\alpha)$ from Remark
\ref{rem-1}. To prove $d\geq 4(\alpha +1)$, we will derive a
contradiction from assuming $d<4(1+\alpha)$.

\
\

{\bf Case 1}: $c_1>0$ and $c_2\geq 0$.

In this case $c(x)\geq 0$, where $c(x)$ is a function defined as
in the proof of Proposition \ref{prop-bdth}. We assume
$2(1+\alpha) \leq d<4(1+\alpha)$ by contradiction. For any
$\lambda>0$, set $E_\lambda=\{x\in \R^2_+:|x|>\frac
1{\sqrt{\lambda}}\}$ and $u_\lambda (x)=u(\lambda
x)+2(1+\alpha)\ln\lambda$.  Then $u_\lambda (x)$ satisfies \beq
\label{eq-4} \left\{
\begin{array}{rcl}
-{\Delta} u_\lambda (x) &=& |x|^{2\alpha}e^{u_{\lambda}}, \qquad
\text{in}
\quad E_\lambda\\
\frac {\partial u_\lambda}{\partial t}&=& c(x)
|x|^{\alpha}e^{\frac {u_\lambda}2}, \qquad \text{on}\quad \partial
E_\lambda \cap
\partial \R^{2}_{+}.
\end{array}
\right. \eeq

\noindent Set
\begin{eqnarray*}
v_\lambda (x)&=& v(\lambda x)+2(1+\alpha)\ln\lambda\\
&=&u(\frac {x}{\lambda |x|^2}) +2(\alpha +1)\ln
\frac 1{\lambda |x|^2}\\
\end{eqnarray*}
where $v(x)$ is the Kelvin transformat of $u(x)$, i.e. $
v(x)=u(\frac{x}{|x|^2}) -4(\alpha +1)\ln |x|$. So, $v_\lambda (x)$
is also a solution of (\ref{eq-4}).

Set $w_\lambda=u-v_\lambda$. Since $E_\lambda $ does not contain
the point $x=0$, $w_\lambda$ is smooth in $E_\lambda$, and
$w_\lambda$ satisfies

\beq\label{eq-5} \left \{
\begin{array}{rcll} -{\Delta} w_\lambda
(x)&=& c_1(x)|x|^{2\alpha}w_\lambda, &\qquad \text{in}
\quad E_\lambda\\
\frac {\partial w_\lambda}{\partial t}&=& c(x)
c_2(x)|x|^{\alpha}w_\lambda, &\qquad \text{on}\quad
\partial E_\lambda \cap
\partial \R^{2}_{+}\\
w_\lambda &=& 0, &\qquad \text{on} \quad \partial E_\lambda\cap
\{t>0\}.
\end{array}
\right. \eeq where $c_1(x)=e^{\xi_1(x)}$ and $c_2(x)=\frac 12
e^{\frac {\xi_2(x)}2}$, $\xi_i(i=1,2)$ are two functions between
$u$ and $v_\lambda$.

\ \

\noindent {\bf Claim 1.}  For $\lambda$ large enough, $w_\lambda
(x)\geq 0$ for all $x\in E_\lambda$.

\ \

{\bf Step 1.}  $\exists R_0$ such that for all $x\in \{x\in
\R^2_+, \frac 2{\sqrt{\lambda}} \leq |x|\leq  R_0\}$, we have
$w_\lambda\geq 0$.

For $x\in \{x\in \R^2_+, \frac 2{\sqrt{\lambda}}\leq|x|\leq R_0\}$
with $R_0$ small enough, we have
\begin{eqnarray*}
w_\lambda (x)&=& u(x)-u(\frac x{\lambda
|x|^2})+2(\alpha+1)\ln(\lambda|x|^2)\\
&\geq & o(1)+2(\alpha+1)\ln 4>0.
\end{eqnarray*}

\ \

{\bf Step 2.} $ \exists R_1\leq R_0$ such that for all $x\in
\{x\in \R^2_+,\frac 1{\sqrt{\lambda}}\leq |x|\leq \frac
2{\sqrt{\lambda}}\leq R_1\}$, we have $w_\lambda\geq 0$.

Set $A_\lambda =\{x\in \R^2_+,\frac 1{\sqrt{\lambda}}\leq |x|\leq
\frac 2{\sqrt{\lambda}}\leq R_1\}$ and $g(x)=1-|x|^{\alpha+1}$ and
let $\overline{w}_\lambda (x)=\frac{w_\lambda (x)}{g(x)}$. Then,
by step 1 and (\ref{eq-5}), $\overline{w}_\lambda (x)$ satisfies
\beq\label{eq-6} \left\{
\begin{array}{lcl}
{\Delta} \overline{w}_\lambda (x)+\frac 2g\nabla g\cdot\nabla
\overline{w}_\lambda (x)+(c_1(x)|x|^{2\alpha}+\frac {{\Delta}
g}{g})\overline{w}_\lambda (x)=0, & &\text{in} \quad
A_\lambda\\
\frac {\partial \overline{w}_\lambda (x)}{\partial
t}=c(x)c_2(x)|x|^\alpha\overline{w}_\lambda (x), & &\text{on}
\quad
\partial A_\lambda\cap \{t=0\}\\
\overline{w}_\lambda\geq 0, & &\text{on} \quad
\partial A_\lambda\cap \{t>0\}\\
\end{array}
\right. \eeq

Since $v_\lambda \leq \max_{\R^2_+}{u}$ in $\overline
{E}_\lambda$,  there exists some positive constant $C_0$ such
that $c_1(x)\leq C_0$. By a direct computation,

\begin{eqnarray*}
c_1(x)|x|^{2\alpha}+\frac {{\Delta} g}g &\leq &
g^{-1}(-(\alpha+1)^2|x|^{\alpha-1}+C_0|x|^{2\alpha}(1-|x|^{\alpha+1}))\\
&\leq & g^{-1}|x|^{\alpha -1}(-(\alpha+1)^2+C_0|x|^{\alpha+1})<0,
\end{eqnarray*}
if $|x|<\{\frac{(\alpha+1)^2}{C_0}\}^{\frac 1{\alpha+1}}$.
Therefore, we choose
$R_1<\min\{\{\frac{(\alpha+1)^2}{C_0}\}^{\frac 1{\alpha+1}},1\}$
small enough. Then, from (\ref{eq-6}), it follows from the maximum
principle and the Hopf Lemma that $w_\lambda \geq 0$ in
$A_\lambda$. Here we have used the fact that $c(x)\geq 0$.

\ \

{\bf Step 3.} $ \exists R_2\leq R_1$ such that for
$\sqrt{\lambda}\geq \frac{ 1}{R_2}$, we have $w_\lambda\geq 0$ for
all $x\in \{x\in \R^2_+, |x|>R_0\}$.

For $x\in \{x\in R^2_+, |x|>R_0\}$ and $d<4\alpha+4$,  as
$|x|\rightarrow \infty$ we have
$$
\lim_{|x|\rightarrow \infty}\frac{u(x)+4(\alpha +1)\ln |x|}{\ln
|x|}=-d+4(\alpha+1)>0.
$$
Then there exists some constant $C>0$ such that
$$
u(x)+4(\alpha +1)\ln |x|>-C,\qquad \text{for}\quad |x|>R_0.
$$
Therefore, for $\lambda$ large enough we have
\begin{eqnarray*}
w_\lambda(x)&=& u(x)+4(\alpha +1)\ln |x|-u(\frac x{\lambda
|x|^2})+2(\alpha+1)\ln \lambda \\
&\geq & -C-\max_{\R^2_+}u+2(\alpha+1)\ln \lambda\geq 0.
\end{eqnarray*}

\noindent Thus we finish the proof of Claim 1.

\ \

 Now we define
$$
\lambda_0=\inf\{\lambda>0|w_\mu (x)\geq 0 \quad \text{in } E_\mu
 \text{ for all }\mu\geq \lambda\}.
$$
\ \

{\bf Claim 2.} $\lambda_0>0$

Assume by contradiction that $\lambda_0=0$, that is, for all
$\lambda >0$, we have $w_\lambda (x)\geq 0$ in $E_\lambda$. Then,
we have for all $x\in \R^2_+$

\begin{equation*}
\left\{
\begin{array}{lcl}
 w_{\frac {1}{|x|^2}}(x) = 0, \\
 w_{\frac {1}{|x|^2}}(rx) \geq 0, \qquad \forall 0<r<1.
\end{array}
\right.
\end{equation*}
Since
$$
w_\lambda(x)=u(x)-u(\frac x{\lambda |x|^2})+2(\alpha+1)\ln
(\lambda |x|^2),
$$
by a direct computation, we have
\begin{equation}\label{111}
w_{\frac 1{|x|^2}}(rx)=u(rx)-u(\frac xr)+4(\alpha+1)\ln r.
\end{equation}
In (\ref{111}), taking firstly $|x|=r$ and then let $r\rightarrow
0^+ $, we get $w_{\frac 1{|x|^2}}(rx)\rightarrow -\infty$. Thus we
get a contradiction with $w_{\frac 1{|x|^2}}(rx)\geq 0 $ for all
$0<r<1$ and all $x\in \R^2_+$.

\ \

{\bf Claim 3.} $w_{\lambda_0}(x)=0, \forall x \in \R^2_+$.

Assume by contradiction $w_{\lambda_0}\geq 0 $ for all $x\in
\R^2_+$. Then from (\ref{eq-5}) we obtain firstly

\begin{equation}\label{eq-8}\left \{
\begin{array}{lcl} {\Delta} w_{\lambda_0}
(x)&\leq & 0, \qquad \text{in}
\quad E_{\lambda_0}\\
\frac {\partial w_{\lambda_0}}{\partial t}&\geq & 0, \qquad
\text{on}\quad
\partial E_{\lambda_0} \cap
\partial \R^{2}_{+}\\
w_{\lambda_0} &=& 0, \qquad \text{on} \quad \partial
E_{\lambda_0}\cap \{t>0\}.
\end{array}
\right. \end{equation} Then we use the strong maximum
principle and the Hopf Lemma to obtain

\beq\label{eq-7} \left \{
\begin{array}{lcl}
w_{\lambda_0} (x)& > & 0, \qquad \text{in}
\quad E_{\lambda_0}\\
\frac {\partial w_{\lambda_0}}{\partial \nu}& > & 0, \qquad
\text{on}\quad
\partial E_{\lambda_0} \cap
\{t>0\}
\end{array}
\right. \eeq where $\nu$ denotes the outward unit normal of the
surface $\partial B_{\sqrt{\frac 1{\lambda_0}}}(0)\cap \{t>0\}$.

Next note that by the definition of $\lambda_0$, we can assume that
there exists a sequence $\lambda_k\rightarrow \lambda_0$ with
$\lambda_k<\lambda_0$ such that $$ \inf_{E_{\lambda_k}}
w_{\lambda_k}<0.$$

If we can prove that \beq\label{117} w_{\lambda_0}(x)\geq C\qquad
\text{for}\quad x\in \overline{E}_{\frac {\lambda_0}2} \eeq for
some constant $C=C(\lambda_0)>0$, then from the continuity of $u$
at $x=0$ we get
$$
w_{\lambda_k}\geq \frac {C}{2}, \qquad\forall x\in
\overline{E}_{\frac {\lambda_0}2}.
$$
for $k$ large enough. It follows that there exists $x_k=(s_k,
t_k)\in \overline{E}_{\lambda_k}\setminus E_{\frac {\lambda_0}2}$
such that
$$
w_{\lambda_k}(x_k)=\inf_{E_{\lambda_k}}w_{\lambda_k}<0.
$$

It is clear that $\sqrt{\frac 1{\lambda_k}}<|x_k|<\sqrt{\frac
2{\lambda_0}}$ and , due to the boundary condition, $t_k>0$. Hence
$\nabla w_{\lambda_k}(x_k)=0$. After passing to a subsequence (still
denoted as $x_k$) $x_k\rightarrow x_0=(s_0, t_0)$, it follows that
\begin{equation}\label{112}
w_{\lambda_0}(x_0)=0, \qquad \nabla w_{\lambda_0}(x_0)=0.
\end{equation}
By (\ref{eq-7}) we have $t_0=0$ and $|s_0|=\sqrt{\frac
1{\lambda_0}}$.

However, we would like to show \beq\label{116}\frac {\partial
w_{\lambda_0}(x_0)}{\partial \nu}>0,
 \qquad \text{for} \quad x_0=(s_0, 0), |s_0|=\sqrt{\frac 1{\lambda_0}}\eeq
if $w_{\lambda_0}(x)$ satisfies (\ref{eq-8}). Here $\nu$ denotes
the outward unit normal of the surface $\partial B_{\sqrt{\frac
1{\lambda_0}}}(0)\cap \{t\geq 0\}$.

Therefore from (\ref{112}) and (\ref{116}) we get a contradiction.
Thus to prove  Claim 3,  it suffices to show (\ref{117})
and (\ref{116}).

 {\bf Proof of (\ref{117})} First, for $x\in
\overline{E}_{\frac {\lambda_0}2}$, we have
\begin{eqnarray*}
v_{\lambda_0}&=& u(\frac x{\lambda_0 |x|^2})+2(\alpha+1)\ln \frac
1{\lambda_0 |x|^2}\\
&\leq & \max_{\R^2_+}u+2(\alpha+1)\ln 2\leq C.
\end{eqnarray*}

Notice that $\min_{\partial E_{\frac {\lambda_0} 2}\cap
\{t>0\}}w_{\lambda_0}\geq \varepsilon$ for some $0<\varepsilon
<1$. Without loss of generality, we assume $\lambda_0=2$. For
$0<r<1$, we introduce an auxiliary function
$$
\varphi(x)=\frac{\varepsilon \mu}{2(c+1)}-\frac {\log|x|}{
\log\sqrt{\frac 1r}}\cdot \varepsilon +\frac {\varepsilon
(1-\mu)(t-\frac 1{\sqrt {r}})}2 (\sqrt{\frac 1r})^\alpha
$$
when $\alpha \geq 0$. Here $c=\max\{c_1,c_2\}$; $0<\mu<1$ will be
chosen later. When $-1<\alpha <0$, we use instead the auxiliary
function
$$
\varphi(x)=\frac{\varepsilon \mu}{2(c+1)}-\frac {\log|x|}{
\log\sqrt{\frac 1r}}\cdot \varepsilon +\frac {\varepsilon
(1-\mu)t}2.
$$
We shall only present the details for the case $\alpha \geq 0$ as the case $-1<\alpha
<0$ can be treated in a similar way. Let
$P(x)=w_{\lambda_0}(x)-\varphi(x)$. Then we get
\begin{equation*}
\left\{
\begin{array}{lcl}
{\Delta} P(x)={\Delta} w_{\lambda_0}(x)\leq 0, & &\qquad
\text{in}\quad E_1\setminus E_r\\
\frac {\partial P(x)}{\partial t}=c(x)c_2(x)|x|^\alpha
w_{\lambda_0}-\frac {\varepsilon(1-\mu)}2(\sqrt{\frac 1r})^\alpha,
&&\qquad \text{on}\quad \partial (E_1\setminus E_r)\cap\{t=0\}.
\end{array}
\right.
\end{equation*}
We will show \beq\label{115} P(x)\geq 0, \quad x\in E_1\setminus
E_r. \eeq We prove it by contradiction. If (\ref{115}) does not
hold, there exists some $x_0=(s_0,t_0)$ such that
$$
P(x_0)=\min_{\overline{E}_1\setminus E_r}P(x)<0.
$$
Since we have $P(x)\geq 0$ on $\partial E_1\cap \{t>0\}$ and
$P(x)>w_{\lambda_0}(x)\geq 0$ on $\partial E_r\cap \{t>0\}$, then
it follows from the maximum principle that $t_0=0$ and
$1<|s_0|<\sqrt{\frac 1r}$ and $\frac {\partial P(x)}{\partial
t}|_{x_0}\geq 0$.

In virtue of $P(x_0)<0$ and $v_{\lambda_0}(x_0)<C_1$, we have
$c_2(x_0)<C_0$ for some constant $C_0>0$ and moreover
\begin{equation}\label{113}
w_{\lambda_0}(x_0)<\varphi(x_0)<\frac {\varepsilon \mu}{2(c+1)}.
\end{equation}

On the other hand, in virtue of $\frac {\partial P(x)}{\partial
t}|_{x_0}\geq 0$ we have
\begin{eqnarray*}
0 &\leq & (c+1)c_2(x_0)|x_0|^\alpha w_{\lambda_0}(x_0)-\frac
{\varepsilon(1-\mu)}2(\sqrt{\frac 1r})^\alpha\\
&\leq &
\{\sqrt{\frac 1r}\}^\alpha (C_0(c+1)w_{\lambda_0}(x_0)-\frac
{\varepsilon(1-\mu)}2)
\end{eqnarray*}
Hence
\begin{equation}\label{114}
w_{\lambda_0}(x_0)\geq \frac {\varepsilon(1-\mu)}{2 C_0(c+1)}
\end{equation}

From (\ref{113}) and (\ref{114}), we have
$$
\mu>\frac 1{1+C_0}.
$$
If we choose $\mu$ such that $0<\mu<\frac 1{1+C_0}$ from the
beginning we reach a contradiction.

Since $P(x)\geq 0$, we then let $r\rightarrow 0$ and have proved
(\ref{117}) with $C=\frac \varepsilon {2(1+c)(1+C_0)}$.

\ \

{\bf Proof of (\ref{116})} Without loss of generality, we assume
$\lambda_0=1$ and $s_0=1$. Set $\Omega=\{x=(s, t)|
1<s^2+t^2<4,\quad s>0,\quad 0<t<\frac 14\}$. Let
$$
h(x)=\varepsilon(s-1)(t+\mu),
$$
and
$$
g(x)=h(x)-h(\frac x{|x|^2})
$$
where $0<\varepsilon, \mu<1$ are chosen later. A direct
computation yields ${\Delta} g(x)=0$ for $ x\in \Omega$. Now
consider
$$
f(x)=w_{\lambda_0}(x)-g(x).
$$
Then we have
\begin{equation*}
\left\{
\begin{array}{lcl}
{\Delta} f(x)\geq 0,  && \text{in}\quad\Omega\\
\frac{\partial f(x)}{\partial t}=c(x)c_2(x)|x|^\alpha
w_{\lambda_0}(x)-\frac{\partial g(x)}{\partial t}, &&
\text{on}\quad \partial \Omega\cap \{t=0\}.
\end{array}
\right.
\end{equation*}

Next we want to show
$$
f(x)\geq 0,\qquad \forall x\in \Omega,
$$
for suitably chosen $\varepsilon$ and $\mu$.

In fact,  we argue by the contradiction and assume that there
exists some $x_1=(s_1, t_1)\in \overline{\Omega}$ such that
$$
f(x_1)=\min_{\overline{\Omega}}f(x)<0.
$$
Since $f(x)=0$ on $\partial \Omega\cap \partial E_1$ and $f(x)\geq
0 $ on $\partial \Omega\cap \{\partial E_{\frac 14}\cup
\{t=\frac14\}\}$,  we can use the maximum principle to obtain
$t_1=0$, $1<s_1<2$ and $\frac {\partial f(x_1)}{\partial t}\geq 0$
on $\partial \Omega \cap\{t=0\}$.

A simple calculation yields
$$
\frac {\partial g(x_1)}{\partial t}=\varepsilon
(s_1-1)(1+s^{-3}_1).
$$
In virtue of $\frac {\partial f(x_1)}{\partial t}\geq 0$ on
$\partial \Omega \cap\{t=0\}$, we obtain
$$
cc_2(x_1)|s_1|^{\alpha}w_{\lambda_0}(x_1)\geq \varepsilon
(s_1-1)(1+s^{-3}_1).
$$
Hence, we get
$$
2^{\alpha}cc_2(x_1)w_{\lambda_0}(x_1)\geq \varepsilon
(s_1-1)(1+s^{-3}_1).
$$
for $\alpha \geq 0$. And
$$
cc_2(x_1)w_{\lambda_0}(x_1)\geq \varepsilon (s_1-1)(1+s^{-3}_1)
$$
for $-1<\alpha <0$. Here $c=\max\{c_1,c_2\}$. On the other hand,
we have
$$
w_{\lambda_0}(x_1)<f(x_1)=\varepsilon \mu(s_1-1)(1+\frac 1{s_1}).
$$
Therefore, if $\alpha \geq 0$ we have
$$
2^{\alpha}(1+cc_2(x_1))\mu\geq \frac{1+s^{-3}_1}{1+s^{-1}_1}>\frac
34,
$$
and if $-1<\alpha <0$, we have
$$
(1+cc_2(x_1))\mu\geq \frac{1+s^{-3}_1}{1+s^{-1}_1}>\frac 34,
$$
If we choose $\mu$ such that $0<\mu<\frac
3{2^{a+2}(1+c\sup_{\R^2_+}c_2(x))}$ for $a=\max\{\alpha,0\}$ from
the beginning we reach a contradiction. Thus we have proved that
$f(x)\geq 0$ for $x\in \Omega$. Since $f(x_0)=0$, i.e. $x_0$ is
minimum point of $f(x)$ in $\overline{\Omega}$, it follows from the
Hopf Lemma that
$$
\frac {\partial f(x_0)}{\partial \nu}\geq 0.
$$
A direct calculation shows that
$$
\frac {\partial w_{\lambda_0}(x_0)}{\partial \nu}=\frac {\partial
f(x_0)}{\partial \nu}+\frac {\partial g(x_0)}{\partial \nu}\geq
\frac {\partial g(x_0)}{\partial \nu}=2\varepsilon \mu>0.
$$
We finish the proof of (\ref{116}).

\
\

In claim 3, $w_{\lambda_0}(x)=0$ implies that
\begin{equation}\label{11}
u(x)=u(\frac {x}{\lambda_0 |x|^2}) +2(\alpha +1)\ln \frac
1{\lambda_0 |x|^2}.
\end{equation}
Hence it follows from (\ref{11}) that $d=4(1+\alpha)$. This
contradicts our assumption $d<4(1+\alpha)$. Thus we proved $d\geq
4(1+\alpha)$. From Proposition \ref{prop-dth2} we know
$d=4(1+\alpha)$.

\
\

{\bf Case 2. $c_1>0$ and $c_2<0$.}

In this case, we will follow the argument of the case 1. The main
difference between the case $c_2\geq 0$ and $c_2<0$, in view of the
maximum principle and the Hopf Lemma, is to show step 2
in the proof of Claim 1. Actually we can prove this step in the
case $c_2<0$ by using a suitable test function. This will become
evident from the rest of the argument.

\ \

{\bf Step 2 of Claim 1}:  $ \exists R_1\leq R_0$ such that for all
$x\in A_\lambda=\{x\in \R^2_+,\frac 1{\sqrt{\lambda}}\leq |x|\leq
\frac 2{\sqrt{\lambda}}\leq R_1\}$, we have $w_\lambda\geq 0$.

Let $x=(s,t)$ and $z=x+(0, +\frac \mu{\sqrt{\lambda}})$, where
$\mu$ is a positive number that will be determined later. Set
$g(x)=1-|z|^{\alpha+1}$ and $\overline{w}_\lambda
(x)=\frac{w_\lambda (x)}{g(x)}$. Then, by step 1 and (\ref{eq-5}),
$\overline{w}_\lambda (x)$ satisfies \beq\label{eq-9} \left\{
\begin{array}{lcl}
{\Delta} \overline{w}_\lambda (x)+\frac 2g\nabla g\cdot\nabla
\overline{w}_\lambda (x)+(c_1(x)|x|^{2\alpha}+\frac {{\Delta}
g}{g})\overline{w}_\lambda (x)=0, \quad \text{in} \quad
A_\lambda\\
\frac {\partial \overline{w}_\lambda (x)}{\partial
t}=(c_1c_2(x)|x|^\alpha-\frac 1g \frac {\partial g}{\partial
t})\overline{w}_\lambda (x), \quad \text{on} \quad
\partial A_\lambda\cap \{t=0\}\cap \{s>0\}\\
\frac {\partial \overline{w}_\lambda (x)}{\partial
t}=(c_2c_2(x)|x|^\alpha-\frac 1g \frac {\partial g}{\partial
t})\overline{w}_\lambda (x), \quad \text{on} \quad
\partial A_\lambda\cap \{t=0\}\cap\{s<0\}\\
\overline{w}_\lambda\geq 0,\quad \text{on} \quad
\partial A_\lambda\cap \{t>0\}\\
\end{array}
\right. \eeq

Since $v_\lambda \leq \max_{\R^2_+}{u}$ in $\overline
{E}_\lambda$, then there exists some positive constant $C_0$ such
that $0<c_1(x), c_2(x)\leq C_0$. Since $x\in A_\lambda=\{x\in
\R^2_+,\frac 1{\sqrt{\lambda}}\leq |x|\leq \frac
2{\sqrt{\lambda}}\leq R_1\}$, we have $|x|\sim |z|\sim |\frac
1{\sqrt{\lambda}}|$.  Then by a direct computation, we obtain

\begin{eqnarray*}
c_1(x)|x|^{2\alpha}+\frac {{\Delta} g}g &\leq &
g^{-1}(-(\alpha+1)^2|z|^{\alpha-1}+C_0|x|^{2\alpha}(1-|z|^{\alpha+1}))<0,
\end{eqnarray*}
if $\lambda $ is large enough. Similarly,  we have

\begin{eqnarray*}
&& c_2c_2(x)|x|^{2\alpha}-\frac 1g\frac {\partial g}{\partial t}
\\
&\geq & g^{-1}((\alpha+1)|z|^{\alpha-1}\frac \mu{\sqrt{\lambda}}
+c_2C_0|x|^\alpha(1-|z|^{\alpha+1}))\\
&\geq & g^{-1}((\alpha+1)C\mu|\frac
1{\sqrt{\lambda}}|^\alpha+c_2C_0|\frac 1{\sqrt{\lambda}}|^\alpha)
>0,
\end{eqnarray*}
on $\partial A_\lambda\cap \{t=0\}\cap\{s<0\}$ for sufficiently large $\mu$. It is obvious that $ c_1c_2(x)|x|^{2\alpha}-\frac 1g\frac
{\partial g}{\partial t}>0$ on $\partial A_\lambda\cap
\{t=0\}\cap\{s>0\}$ since $c_1>0$. Then, from (\ref{eq-9}), we can
again use the maximum principle and the Hopf Lemma to obtain $w_\lambda
\geq 0$ in $A_\lambda$.

The proof of Claim 3 requires some simple modifications
when we use the maximum principle and the Hopf Lemma. But these can be
carried out just by  changing test functions as in the previous argument.
Here we omit the details. Thus we complete the proof of the Theorem.
\end{proof}

\begin{rem}
Actually the spherical symmetry (\ref{11}) is inherited by the
solution of (\ref{eq-1})-(\ref{con-1}). From the proof of
Proposition \ref{prop-dth3}, it is sufficient to establish Step 3
when $d=4(1+\alpha)$. But this can be done with the help of the
asymptotic estimate (\ref{asy-est4}).
\end{rem}

\
\

\section {Proof of Main Theorems}
In this section we prove our main theorems. Theorem \ref{thm1} can
be obtained directly from Proposition \ref{prop-dth3}, since we
can show that the solution $u$ to (\ref{eq-1})-(\ref{con-1}) has a
removable singularity at $z=\infty$ by using the Kelvin
transformation as in many conformal problems. To prove Theorem
\ref{theo-metr}, we follow closely the argument in \cite{T2}. The
crucial step is to construct a projective connection on $S^2$ by
using the conformal metric on $\overline{\R^2_{+}}\cup\{\infty\}$
with constant curvature $1$ and constant geodesic  curvature
$c(x)$ on the boundary.

First, we prove Theorem \ref{thm1}:

{\bf Proof of Theorem \ref{thm1}} To prove Theorem \ref{thm1}, it
 suffices to show that any solution of (\ref{eq-1})
-(\ref{con-1}) determines a metric as in Theorem \ref{thm1}. For
this point, we first prove that the metric
$ds^2=|x|^{2\alpha}e^{u(x)}|dz|^2$,  $u$ being a solution of
(\ref{eq-1}) -(\ref{con-1}), has two conical singularities at $0$
and $\infty$ with the same order. The existence of this
metric is shown in Theorem \ref{theo-metr}.

Let $v$ be the Kelvin transformation of $u$. If $u$ is a solution
of (\ref{eq-1}) -(\ref{con-1}), then $v\in C^2(\R^{2}_{+})\cap
C^1(\overline{\R^{2}_{+}}\setminus\{0\})$ and satisfies \beq
\left\{
\begin{array}{rcl}
-{\Delta} v &=& |x|^{2\alpha}e^v, \qquad \text{in}
\quad \R^{2}_{+},\\
\frac {\partial v}{\partial t}&=& c(x) |x|^{\alpha}e^{\frac v2},
\qquad \text{on}\quad \partial \R^{2}_{+}\setminus\{0\}.
\end{array}
\right. \eeq To prove the result, we first show that $v$ is
continuous at $x=0$, that is the singularity $z=0$ of $v$ is
removable. Applying the asymptotic estimate (\ref{asy-est4}) we
have
\begin{eqnarray*}
v(x)&=& u(\frac{x}{|x|^2}) -4(\alpha +1)\ln |x|\\
&=&(d-4(\alpha +1))\ln{|x|}+O(1) \qquad \text{for}\quad |x| \quad
\text{near}\quad 0.
\end{eqnarray*}
Since $d =4(1+\alpha)$, we get that $v$ is bounded near $0$. Thus, by standard elliptic regularity, we conclude that $v$ is a
$C^2(\R^{2}_{+})\cap C^1(\overline{\R^{2}_{+}})$ solution of
(\ref{eq-1}) when $\alpha \geq 0$. While, for $\alpha \in (-1,0)$,
$v$ is smooth away from the origin and $v\in W^{2,p}$ for
$1<p<-\frac 1{\alpha}$ near the origin. In particular, in any case, $v$ is continuous at the origin.

Next note that $ds^2=e^{\widetilde{u}}dx^2$ for
$\widetilde{u}=u(x)+2\alpha\log|x|$, where $u$ is a solution of
(\ref{eq-1}) -(\ref{con-1}). So the metric $ds^2$ has a conical
singularity at $z=0$ with order $\alpha$. Let
$\widetilde{v}(x)=\widetilde{u}(\frac x{|x|^2})-4\log|x|$ be
Kelvin transformation of $\widetilde{u}$. Then we obtain near $z=0$
\begin{eqnarray*}
\widetilde{v}(x)&=& u(\frac x{|x|^2})-2\alpha\log|x|-4\log|x|\\
&=& 2\alpha\log|x|+v(x)\\
\end{eqnarray*}
since $v(x)$ is continuous function at $z=0$. By the definition of a
conical singularity, we get that the metric
$ds^2=e^{\widetilde{u}}dx^2$  has a conical singularity at
$z=\infty$ with the same order as at $z=0$.\qed

\
\

\begin{lm}\label{lm-Proj}Let $u$ be a solution of
(\ref{eq-1})-(\ref{con-1}), and $ds^2=e^{\tilde{u}}|dz|^2$, where
$\tilde{u}=u+2\alpha \ln |z|$ . Define
$$ \eta(z)=(\frac{\partial^2 \tilde{u}}{\partial z^2}-\frac
12(\frac{\partial\tilde{u}}{\partial z})^2)|dz|^2 .$$ Then
$\eta(z)$ can be extended to a projective connection on
$\S^2=\C\cup \infty$, still denoted by $\eta(z)$, that is compatible
with the divisor {\bf A}$=\alpha\cdot 0+\alpha\cdot\infty$.
\end{lm}

\begin{proof} First, from the assumption, we know that  $\tilde{u}$ satisfies \beq \label{4.1} \left\{
\begin{array}{rcll}
-{\Delta} \widetilde{u} &=& e^{\widetilde{u}}, &\qquad \text{in}
\quad \R^{2}_{+},\\
\frac {\partial \widetilde{u}}{\partial t}&=& c(x) e^{\frac
{\widetilde{u}}{2}}, &\qquad \text{on}\quad
\partial \R^{2}_{+}\setminus\{0\},
\end{array}
\right. \eeq with the energy conditions

\beq \label{4.2} \int_{\R^{2}_{+}}e^{\widetilde{u}}dx<\infty ,
\eeq

\beq \label{4.3} \int_{\partial \R^{2}_{+}}e^{\frac
{\widetilde{u}}{2}}dt<\infty. \eeq

\

Let $f(z)=\frac{\partial^2 \tilde{u}}{\partial z^2}-\frac
12(\frac{\partial\tilde{u}}{\partial z})^2$, then from
(\ref{4.1}), $f(z)$ is holomorphic on $\R^2_{+}$ and ${\bf
\text{Im}}f =\frac 12(\frac 12 \frac{\partial \tilde{u}}{\partial
s}\frac{\partial \tilde{u}}{\partial t}-\frac
{\partial^2\tilde{u}}{\partial s
\partial t})$. On the other hand, since on $\partial
\R^2_+\setminus\{0\}$, $\frac {\partial \widetilde{u}}{\partial
t}= c(z) e^{\frac {\widetilde{u}}{2}}$, we have $\frac
{\partial^2\tilde{u}}{\partial s \partial t}=\frac {c(z)}2
e^{\frac {\tilde{u}}2}\frac{\partial \tilde{u}}{\partial s}=\frac
12 \frac{\partial \tilde{u}}{\partial s}\frac{\partial u}{\partial
t}$. This implies $f(z)$ is real on $\partial
\R^2_+\setminus\{0\}$, and we may extend $f(z)$ to a holomorphic
function on $\C\setminus\{0\}$ by
$f(z)=\overline{f(\overline{z})}$ for $z\in \R^2_-$. Thus we
extend $\eta$ to $\C$ such that $\eta$ is holomorphic on
$\C-\{0\}$.

\
\

Next we show that $\eta(z)$ is a projective connection on
$\C\cup\infty$. Let $(V,w)$ and
$(U,z)$ be  coordinate charts with $U\cap V\neq \emptyset $. If
$U\cap V \subset \overline{\R^2_+}\cup \{\infty\}$, then by
following the argument in \cite{T2} and by using the fact that
$ds^2=e^{\tilde{u}}|dz|^2$ is a conformal metric on $\R^2_+\cup
\infty$, we have $ds^2=e^{\tilde{u}}|dz|^2=e^v|dw|^2$ with
$v=\tilde{u}+\frac 12\log|\frac{dz}{dw}|$, and consequently we get
\begin{eqnarray}\label{4.4}
\eta(w)&=&(\frac {\partial^2(\tilde{u}+\frac
12\log|\frac{dz}{dw}|)}{\partial w^2}-\frac 12(\frac
{\partial(\tilde{u}+\frac 12\log|\frac{dz}{dw}|)}{\partial
w})^2)|dw|^2 \nonumber\\
& =& \eta(z)+\{z,w\}|dw|^2.
\end{eqnarray}
If  $U\cap V \subset \R^2_-$,  since $\overline
{\bar{z}_{\bar{w}}}=z_w$, we get from (\ref{4.4})
\begin{eqnarray*}
\eta(w)&=&\overline{\eta(\bar{w})}=\overline{\eta(\bar{z})+\{\bar{z},\bar{w}\}d\bar{w}^2}\\
&=&\eta(z)+\{z,w\}|dw|^2.
\end{eqnarray*}
So, in any case,  $\eta(w)=\eta(z)+\{z,w\}dw^2$ when
$U\cap V\neq \emptyset$. This means that $\eta$ is a
projective connection on $S^2=\C\cup\infty$.

\
\

Next, we want to show that $\eta$ has a regular singularity at $0$
and at $\infty$ of weight $\rho=-\frac 12\alpha(\alpha+2)$. We
prove this statement only at the singular point $0$, since  the same
argument applies at $\infty$ by using the Kelvin
transformation. Since $0$ is a conical singular point of the
metric $ds^2=e^{\tilde{u}}dz^2$ on $\R^2_+\cup\{\infty\}$, we set
$\tilde{u}=u(x)+2\alpha\log|x|$ in $B_r(0)\cap \overline{\R^2_+}$,
where $u(x)$ is a continuous  solution of
(\ref{eq-1})-(\ref{con-1}).

First, we consider the case $\alpha\geq 0$. In this case, since
$u(x)$ is a continuous solution of (\ref{eq-1})-(\ref{con-1}),  $u$ is
of class $C^2$ in $\overline{\R^2_+}$ by  classical
elliptic regularity theory. Then we have
$$
\frac {\partial^2\widetilde{u}}{\partial z^2}-\frac
12(\frac{\partial \widetilde{u}}{\partial z})^2=\frac
{\partial^2u}{\partial z^2}-\frac 12(\frac{\partial u}{\partial
z})^2-\frac \alpha z\frac {\partial u}{\partial z}-\frac
{\alpha(\alpha+2)}{2z^2}.
$$ in $\overline{\R^2_+}\setminus\{0\}$. Hence we obtain
$$
\eta(z)=(-\frac {\alpha(\alpha+2)}{2z^2}-\frac
{\alpha}z\frac{\partial u(z)}{\partial z}+\phi(z))dz^2,\quad
\text{ for } z\in \overline{\R^2_+}\setminus\{0\},
$$
$$
\eta(z)=(-\frac {\alpha(\alpha+2)}{2z^2}-\frac
{\alpha}z\overline{\frac{\partial u(\bar{z})}{\partial
\bar{z}}}+\overline{\phi(\bar{z})})dz^2,\quad \text{ for } z\in
\R^2_-,$$ where $\phi(z)=\frac {\partial^2u}{\partial z^2}-\frac
12(\frac{\partial u}{\partial z})^2$ for $z\in
\overline{\R^2_+}\setminus\{0\}$.  This proves that $\eta(z)$ has
a regular singularity of weight $\rho=-\frac 12\alpha(\alpha+2)$
at  $z=0$ in this case.

\
\

When $-1<\alpha<0$, $u$ might not to be $C^2$ and the computation
above might not work. However, we may take a method used in
\cite{T2} to lift the metric to a local branched cover: We set
$z=w^m (m\in \N)$, then the metric can be lifted in the $w-plane$:
$ds'^2=e^{u'}dw^2$ with $u'=\widetilde{u}+2\log|\frac
{dz}{dw}|=u+2(m(\alpha +1)-1)\log|w|+2\log m$, when $z$ is in the
upper half plane.   Therefore $ds'^2$ has a conical singularity at
$w=0$ of order $\alpha'=m(\alpha+1)-1$. Since equation
(\ref{4.1}) is invariant under  conformal transformations, $u'$
satisfies (\ref{4.1}) in terms of $w$. Now choosing $m$ large
enough, we have $\alpha'>0$. Then we can use the same argument as in
\cite{T2} and the extension technique above to get
$$
\eta(z)=(-\frac {\alpha(\alpha+2)}{2z^2}+\frac {\sigma
}{z}+\phi(z))dz^2
$$
where $\phi(z)$ is holomorphic function.
\end{proof}

 \ \

\noindent{\bf Proof of Theorem \ref{theo-metr}.} From Lemma \ref{lm-Proj},
we know that $\eta (z)$ is a projective connection on $S^2=\C\cup
\{\infty\}$ with regular singularities at $z=0$ and $z=\infty$. It
follows from Proposition 2 in  \cite{T2} that
$$\eta(z)=-\frac{\alpha(\alpha +2)}2\cdot \frac {dz^2}{z^2}$$ in
the standard coordinate $z$.

\
\

Setting $h=e^{-\frac {\widetilde{u}}2}$, then we have
\begin{equation}\label{eq-10}
\frac{\partial^2h}{\partial z^2}=\frac {\alpha(\alpha+2)}4\cdot
\frac h{z^2},  \quad \text{ for any } z\in \R^2_+,
\end{equation}
and the boundary condition is \begin{equation}\label{con-5}
\frac{\partial h}{\partial \bar{z}}-\frac{\partial h}{\partial
z}=-\frac {ic(x)}2, \quad\text{ on }\partial \R^2_+\setminus
\{0\}.
\end{equation}
All solutions of (\ref{eq-10}) are of the form
$$
h(z,\bar{z})=f(\bar{z})z^{-\frac{\alpha}2}+g(\bar{z})z^{1+\frac{\alpha}2},
$$
for any $z\in \R^2_+$. Since $h$ is real and analytic, we have
$$
h(z,\bar{z})=a(\bar{z}z)^{-\frac{\alpha}2}+pz^{1+\frac \alpha
2}\bar{z}^{-\frac{\alpha}2}+\bar{p}\bar{z}^{1+\frac \alpha
2}z^{-\frac{\alpha}2}+b(z\bar{z})^{1+\frac{\alpha}2}, \quad\text{
for any }z\in \R^2_+.
$$
Here, $a$ , $b\in \R$ and $p\in \C$.  Since
$\widetilde{u}=u+2\alpha\ln|x|$ near $0$ for some continuous
function $u$, it is clear that $a\neq 0$. Then rewriting
$h(z,\bar{z})$, we have
$$
h=a\cdot( \frac{|1+\bar{\mu} \bar{z}^{\alpha+1}|^2+\nu|z|^{2\alpha
+2}}{|z|^{\alpha}}),
$$
for some parameters $\mu=\frac pa\in \C$ and
$\nu=\frac{ab-p\bar{p}}{a^2}\in \R$. Therefore, a conformal metric
should be
$$
ds^2=\frac{|dz|^2}{h^2}=\frac 1{a^2}\cdot \frac
{|z|^{2\alpha}|dz^2|}{(|1+\bar{\mu}
\bar{z}^{\alpha+1}|^2+\nu|z|^{2\alpha +2})^2}.
$$
Setting $w=\frac 1{\bar{z}}$, we have
$$
ds^2=\frac 1{a^2}\cdot \frac {|w|^{2\alpha}|dw^2|}{(|\bar{\mu}
+w^{\alpha+1}|^2+\nu)^2}.
$$
On the other hand, if we assume $(r,\theta)$ is the polar
coordinate system in $\R^2$, then we have
$$
h(r,\theta)=ar^{-\alpha}+pre^{i\theta(1+\alpha)}+\bar{p}re^{-i\theta(1+\alpha)}+br^{2+\alpha}.
$$
And its boundary condition (\ref{con-5}) can be rewritten as
$$
-\frac{\partial h}{\partial
\theta}(e^{i\theta}+e^{-i\theta})+ir\frac{\partial h}{\partial
r}(e^{i\theta}-e^{-i\theta})=rc(r,\theta),$$ for $\theta=0$ and
$\theta=\pi$. Here $c(r,\theta)=c_1$ if $\theta=0$ and
$c(r,\theta)=c_2$ if $\theta=\pi$. Therefore we obtain by using
the partial derivative $\frac {\partial h}{\partial \theta}$ at
$\theta=0$ and $\theta=\pi$ respectively
$$
2(\alpha+1)(\bar{p}-p)=-ic_1,
$$
and $$ 2(\alpha+1)(\bar{p}e^{-i\alpha\pi}-pe^{i\alpha\pi})=-ic_2.
$$

\
\

Then there are two cases.

\
\

In the first case, $\alpha$ is an integer: When $\alpha=2k, k=0,1,2,\dots$, then $c_1=c_2$. And when $\alpha=2k+1, k=0,1,2,\dots$,
then $c_1=-c_2$.  In this case one can only
determine $\text{Im}\{p\}$, namely  $\text{Im}\{p\}=\frac
{c_1}{4(\alpha+1)}$.
 Now  we set
$\frac{\text{Im}\{p\}}{a}=\frac{c_1\lambda^{\alpha+1}}{\sqrt{2}}$.
Then we have
$$
a=\frac{\sqrt{2}}{4(\alpha+1)\lambda^{\alpha+1}},
$$
and consequently
$$
ds^2=\frac
{8(\alpha+1)^2\lambda^{2(\alpha+1)}|w|^{2\alpha}|dw^2|}{(|w^{\alpha+1}-w_0|^2+\nu)^2},
$$
where $w_0=(x_0,t_0)$ for some real number $x_0$ and
$t_0=\frac{c_1\lambda^{\alpha+1}}{\sqrt{2}}$. Set
$$u=\log \frac
{8(\alpha+1)^2\lambda^{2(\alpha+1)}}{(|w^{\alpha+1}-w_0|^2+\nu)^2}.$$
Then it follows from the definition of the conformal metric that
$u$ is a solution of (\ref{eq-1})-(\ref{con-1}). Hence we have
$\nu=\lambda^{2\alpha+2}$. This implies
$$
ds^2=\frac
{8(\alpha+1)^2\lambda^{2(\alpha+1)}|w|^{2\alpha}|dw^2|}{(|w^{\alpha+1}-w_0|^2+\lambda^{2\alpha+2})^2}.
$$

\
\

In the second case, $\alpha\neq k,k=0,1,2,\dots$. For
any $c_1$ and $c_2$, one can then find a unique complex number $p$. In
this case, we also set
$\frac{\text{Im}\{p\}}{a}=\frac{c_1\lambda^{\alpha+1}}{\sqrt{2}}$.
Then we have
$$
a=\frac{\sqrt{2}}{4(\alpha+1)\lambda^{\alpha+1}},
$$
and consequently we have
$$
ds^2=\frac
{8(\alpha+1)^2\lambda^{2(\alpha+1)}|w|^{2\alpha}|dw^2|}{(|w^{\alpha+1}-w_0|^2+\nu)^2},
$$
where $w_0=(x_0,t_0)$ is a fixed point for
\begin{equation}\label{xx1}
 x_0=\frac
{\lambda^{\alpha+1}(c_1\cos(\pi\alpha)-c_2)}{\sqrt{2}\sin(\pi\alpha)} \quad
\hbox{ and } \quad t_0=\frac{c_1\lambda^{\alpha+1}}{\sqrt{2}}.\end{equation}
 Then as in the
first case, we can get
$$
ds^2=\frac
{8(\alpha+1)^2\lambda^{2(\alpha+1)}|w|^{2\alpha}|dw^2|}{(|w^{\alpha+1}-w_0|^2+\lambda^{2\alpha+2})^2}.
$$
 We complete the proof. \qed

\

Since the domain  $\overline{\R^2_+}\backslash \{0\}$
is simply connected, in this paper we consider $z^{1+\alpha}$ as a well-defined function, even if for non-integer $\alpha$.
In polar coordinates, we have
\[
 e^u=\frac
{8(\alpha+1)^2\lambda^{2(\alpha+1)}}{((r^{1+\alpha}\cos (1+\alpha)\theta-x_0)^2+(r^{1+\alpha}\sin (1+\alpha)\theta-t_0)^2
+\lambda^{2\alpha+2})^2},
\]
where $x_0$ and $t_0$ are given by (\ref{xx1}).

 \ \

\section*{\bf{Acknowledgements}} The work was carried out when the third
author was visiting the Max Planck Institute for Mathematics in
the Sciences. She would like to thank the institute for the
hospitality and the good working conditions.

\end{document}